\documentstyle[12pt]{article}
\newtheorem{tanim}{Definition}
\newtheorem{teo}{Theorem}
\newtheorem{pro}{Proposition}

\begin{document}

\begin{center}
{\Large {\bf Analysis on the 2--Dim Quantum Poincar\'{e} Group 
at Roots of Unity}}
\end{center}

\vspace{1cm}

\noindent
{\bf Hadji Ahmedov }

\vspace{.5cm}

\noindent
{\small Feza G\"ursey Institute, P.O. Box 6, 81220 Cengelk\"{o}y, Istanbul,
Turkey.}

\noindent
{\small E--mail:hagi@gursey.gov.tr } \vspace{2cm}

\noindent
{\bf Abstract:} 2--Dim quantum Poincar\'e Group $E_q(1,1)$ at roots
of unity, its dual $U_q(e(1,1))$ and some of its homogeneous spaces are
introduced. Invariant integrals on $E_q(1,1)$ and its invariant discrete
subgroup $E(1,1\mid p)$ are constructed. $*$--Representations of 
the quantum algebra $U_q(e(1,1))$ constructed in the homogeneous 
space $SO(1,1\mid p)$ are  integrated to the pseudo--unitary 
representations of $E_q(1,1)$ by means of the universal $T$--matrix. 
$U_q(e(1,1))$ is realized on the quantum plane $E_q^{(1,1)}$ and the 
eigenfunctions of the complete set of observables are obtained in 
the angular momentum and momentum basis. The matrix elements of 
the pseudo--unitary irreducible representations are 
given in terms of the cut off q--exponential and  $q$--Bessel functions 
whose properties we also investigate.

\vspace{1cm}

\noindent 
{\bf 1. Introduction}

\vspace{.2cm} \noindent 
Finite dimensional representations of the quantum algebra $U_q(g)$ for real $%
q$ is very similar to the representations of the universal enveloping
algebra $U(g)$ where $g$ is the complex simple Lie algebra \cite
{jim,lus,ros,uen}. Theory of the algebraic quantum group $G_q$ which is the
Hopf algebra of the quantized polynomials on the Lie group $G$ is
essentially the same as that of $G$ too (see \cite{vil} and references
therein ). Matrix elements of the irreducible representations of $G_q$ are
expressed in term of the q--special functions which are the generalization
of the ones related to the Lie group $G$. There also exist an invariant
distance \cite{hagi}, an invariant integral and Peter--Weyl approximation
theorem \cite{woron} on the compact quantum group $G_q$ and its symmetric
spaces.

On the other hand the quantum algebra $U_q(g)$ at roots of unity admits
finite dimensional irreducible representations which have no classical
analogs \cite{cha,kac1,kac,lus1,roc}. Because of the peculiar algebraic
structure of these representations quantum algebras at roots of unity 
have found interesting applications, especially in determining knot 
invariants \cite{res} and in the quantum Hall effect \cite{dayi}. 
Unlike the case of real $q$ theory of the dual space $G_q$ at roots of 
unity is not well established :

(i) what is the structure of the quantum group $G_q$ at roots of unity ?

(ii) what are the q--special functions related to $G_q$ at roots of unity?

(iii) are there invariants (integral, distance ) on $G_q$ at roots of unity?

The problems (i) and (iii) are partially solved for the quantum group 
$SL_q(2,C)$ at roots of unity in \cite{coq,glush}. 
Quantum groups at roots of unity appear to be a natural generalization of
the usual supersymmetry to the fractional one ( FSUSY ) which replaces the 
$Z_2$--grading of the SUSY algebra with a $Z_p$--graded algebra in such 
a way 
that the FSUSY transformation mix elements of all grades (see \cite{dum1}
and references therein ).

The purpose of this paper is to solve the problems (i), (ii) and (iii) for
the 2--dim quantum Poincar\`e group $E_q(1,1)$ at $q^p=1$. This
group is the $Z_p$--graded product of the $p^3$--dimensional invariant $%
E(1,1\mid p)$ and translation $R^2$ subgroups. We define the invariant
integral on $E_q(1,1)$ and demonstrate that all the methods of
representation theory available at generic $q$ can be extended on this group.

In Section 2 we define the quantum Poincar\'e group $E_q(1,1)$ at roots of
unity, its homogeneous spaces $E(1,1\mid p)$, $SO(1,1\mid p)$, $M^{(1,1)}$, 
$E_q^{(1,1)}$ and the dual space $U_q(e(1,1))$.
Section 3 is devoted to the construction of the invariant integral on $%
E_q(1,1)$ and its invariant discrete subgroup $E(1,1\mid p)$.
The irreducible $*$--representation of $U_q(e(1,1))$ constructed in Section
4 are integrated to the pseudo--unitary irreducible representations of $%
E_q(1,1)$ by means of the universal $T$--matrix in Section 5. The matrix
elements of these representations and some of their properties are
investigated in Section 5 also. 
In Section 6 we realize the quantum algebra $U_q(e(1,1))$ on the quantum
plane $E_q^{(1,1)}$ and obtain the eigenfunctions of the complete set of
commuting elements of $U_q(e(1,1))$ in the angular momentum and momentum
basis.

\vspace{1cm} \noindent 
{\bf 2. 2--Dim Quantum Poincar\'{e} Group $E_q(1,1)$ at Roots of Unity}

\vspace{.2cm} \noindent 
Let us start by reviewing the principal facts of the 2--dimensional complex
quantum Euclidean group $E_q(2,C)$ and its dual $U_q(e(2,C))$ \cite{cel}.

The quantum group $E_q(2,C)$ is the Hopf algebra $A(E_q(2,C))$ generated by $%
\eta _{\pm }$ and $\delta ^{\mp 1}$ satisfying the relations 
\begin{equation}
\label{1}\eta _{-}\eta _{+}=q^2\eta _{+}\eta _{-},\ \ \ \eta _{\pm }\delta
=q^2\delta \eta _{\pm }
\end{equation}
and 
\begin{eqnarray}\label{2} 
\Delta (\eta_\pm )=\eta_\pm\otimes 1_A+\delta^{\pm 1} \otimes \eta_\pm , 
\ \ \ \ \ \Delta (\delta )=\delta\otimes \delta , \ \ \ \ \ \ \ \nonumber \\
\varepsilon (\delta^{\pm 1} )=1,\ \ \ \varepsilon (z_\pm)=0, \ \ \ 
S(\delta^{\pm 1} )=\delta ^{\mp 1},\ \ \ S(\eta_\pm)=-\delta ^{\mp1}
\eta_\pm . 
\end{eqnarray}
The quantum algebra $U_q(e(2,C))$ is the Hopf algebra generated by $p_{\pm }$
and $\kappa ^{\pm 1}$ satisfying the relations 
\begin{equation}
\label{5}p_{+}p_{-}=p_{-}p_{+},\ \ \ \ \ \ p_{\pm }\kappa =q^{\mp 1}\kappa
p_{\pm }
\end{equation}
and 
\begin{eqnarray}\label{6}
\Delta (p_\pm )=p_\pm \otimes \kappa +\kappa^{-1}\otimes p_\pm , 
\ \ \ \ \ \Delta (\kappa)=\kappa\otimes \kappa, \ \ \ \ \ \ \ \nonumber \\
\varepsilon (p_\pm )=0 , \ \ \  \varepsilon (\kappa^{\pm 1})=1,
\ \ \ S(p_\pm )=-q^{\pm 1}p_\pm , \ \ \  S(\kappa^{\pm 1})=\kappa^{\mp 1}. 
\end{eqnarray}
The duality pairings between $A(E_q(2,C))$ and $U_q(e(2,C))$ are given by 
\begin{equation}
\label{9}\langle \kappa ^j,\delta ^{j^{\prime }}\rangle =q^{jj^{\prime }},\
\ \ \ \ j,\ j^{\prime }\in Z
\end{equation}
and 
\begin{equation}
\label{10}\langle p_{\pm }^n,\eta _{\pm }^m\rangle =i^nq^{\pm \frac
n2}[n]!\delta _{nm},\ \ \ \ \ \ n,\ m\in N,
\end{equation}
where 
$$
[n]=\frac{q^n-q^{-n}}{q-q^{-1}},\ \ \ \ \ \ [n]!=[1][2]\cdots [n]. 
$$
Since $\Delta $ is a homomorphism (\ref{2}) implies that
\begin{equation}
\label{13}\Delta (\eta _{\pm }^n)=\sum_{m=0}^n\left[ 
\begin{array}{c}
n \\ 
m
\end{array}
\right] _{\pm }\eta _{\pm }^{n-m}\delta ^{\pm m}\otimes \eta _{\pm }^m,
\end{equation}
where 
$$
\left[ 
\begin{array}{c}
n \\ 
m
\end{array}
\right] _{\pm }=q^{\pm m(m-n)}\frac{[n]!}{[n-m]![m]!}. 
$$
The Hopf algebra $A(E_q(2,C))$ has two real forms $A(E_q(2))$ and $%
A(E_q(1,1))$ defined by the involutions 
$$
\delta ^{*}=\delta ^{-1},\ \ \ \ \eta _{\pm }^{*}=\eta _{\mp }\ \ \ {\rm for}%
\ \ \ q\in R 
$$
and 
\begin{equation}
\label{4}\delta ^{*}=\delta ,\ \ \ \ \eta _{\pm }^{*}=\eta _{\pm }\ \ \ {\rm %
for}\ \ \ \mid q\mid =1
\end{equation}
respectively. The 2--dimensional quantum Euclidean group $E_q(2)$ which is
the $*$--Hopf algebra $A(E_q(2))$ was treated in detail in 
\cite{vak,wor,bon}. $A(E_q(1,1))$ is the 2--dimensional quantum 
Poincar\{e group $E_q(1,1)$. The Hopf algebra $U_q(e(2,C))$ has two real 
forms $U_q(e(2))$ and $U_q(e(1,1))$ defined by the involutions 
$$
p_{\pm }^{*}=p_{\mp },\ \ \ \ \kappa ^{*}=\kappa \ \ \ {\rm for}\ \ \ q\in R 
$$
and 
\begin{equation}
\label{8}p_{\pm }^{*}=p_{\pm },\ \ \ \ \kappa ^{*}=\kappa \ \ \ {\rm for}\ \
\ \mid q\mid =1
\end{equation}
respectively.

For future convenience we would like to introduce the convolution product $%
\diamond$. Let $\xi : A\rightarrow V$ be the homomorphic map of a Hopf
algebra $A$ onto a linear space $V$. We set 
$$
\xi\diamond f = (id\otimes \xi )\Delta (f), \ \ \ f\diamond \xi =
(\xi\otimes id )\Delta (f), \ \ \ \xi\diamond \xi = (\xi\otimes \xi )\Delta
. 
$$
Clearly $\xi\diamond f$ and $f\diamond \xi$ belong to $A\otimes V$ and $%
V\otimes A$ respectively; $\xi\diamond \xi$ is homomorphic map of $A\otimes A
$ onto $V\otimes V$.

When $q$ is a root of unity $q^p=1$ (we deal with odd $p$ ) the duality
relations (\ref{9}) and (\ref{10}) become degenerate. To get rid of these
degeneracies we have to redefine the $*$--Hopf algebras $A(E_q(1,1))$ and $%
U_q(e(1,1))$.

To remove the degeneracy in (\ref{9}) we put 
\begin{equation}
\label{15}\delta ^p=1_A
\end{equation}
and 
\begin{equation}
\label{15a}\kappa ^p=1_U.
\end{equation}
Instead of (\ref{9}) we then have 
\begin{equation}
\langle \kappa ^n,\zeta (m)\rangle =\delta _{nm},\ \ \ \ \ n,\ m\in [0,p-1],
\end{equation}
where 
$$
\zeta (m)=\frac 1p\sum_{n=0}^{p-1}q^{-nm}\delta ^n,\ \ \ m\in [0,\ p-1], 
$$
which satisfies the periodicity  property $\zeta (m+pj)=\zeta (m)$, $j\in Z$.

To remove the degeneracy in (\ref{10}) we put 
\begin{equation}
\label{19}\eta _{\pm }^p=0 
\end{equation}
such that new variables $z_{\pm }$ 
\begin{equation}
\label{20}z_{\pm }=\lim _{q^p=1}(-1)^{\frac{p+1}2}\frac{\eta _{\pm }^p}{[p]!}
\end{equation}
are well defined. The above limiting process stems from the work De Concini,
Kac and collaborators, and Lusztig which also appears in two recent
monographs \cite{mon1}, \cite{mon2}, from which it can be traced back to the
original references. The expression (\ref{10}) now reads 
\begin{equation}
\langle p_{\pm }^n,\eta _{\pm }^m\rangle =i^nq^{\pm \frac n2}[n]!\delta
_{nm},\ \ \ n, \ m\in [0,\ p-1] 
\end{equation}
and 
\begin{equation}
\langle P_{\pm }^n,z_{\pm }^m\rangle =i^nn!\delta _{nm},\ \ \ n, \ m\in N, 
\end{equation}
where $P_{\pm }=p_{\pm }^p$. Inspecting (\ref{1}) and (\ref{20}) we conclude
that the new variables $z_{\pm }$ commute with $\eta _{\pm }$ and $\delta $.
By the virtue of (\ref{13}) and (\ref{20}) we obtain 
$$
\Delta (z_{\pm })=z_{\pm }\otimes 1_A+1_A\otimes z_{\pm }+(-1)^{\frac{p+1}%
2}\sum_{n=1}^{p-1}\frac{q^{\pm n^2}}{[p-n]![n]!}\eta _{\pm }^{p-n}\delta
^{\pm n}\otimes \eta _{\pm }^n. 
$$
Moreover, we have 
$$
S(z_{\pm })=-z_{\pm },\ \ \ \varepsilon (z_{\pm })=0,\ \ \ z_{\pm
}^{*}=z_{\pm }.%
$$
At this point we would like to introduce the short hand notation 
$$
\Delta (z)=Z+B,%
$$
where $z=(z_{+},\ z_{-})$, $Z=(Z_{+},\ Z_{-})$, $B=(B_{+},\ B_{-})$ and 
$$
Z_{\pm }=z_{\pm }\otimes 1_A+ 1_A\otimes z_{\pm },\ \ \ \ B_{\pm }=(-1)^{%
\frac{p+1}2} \sum_{n=1}^{p-1} \frac{q^{\pm n^2}}{[p-n]![n]!}\eta _{\pm
}^{p-n}\delta ^{\pm n}\otimes \eta _{\pm }^n.%
$$
Since $B_{\pm }^2=0$ for any function $f$ from the space $C^\infty (R^2)$ of
all infinitely differentiable functions on $R^2$ we have 
\begin{equation}
\label{26}\Delta (f(z))=f(Z)+\frac{df}{dz_{+}}\mid _{z=Z}B_{+}+ \frac{df}{%
dz_{-}}\mid_{z=Z}B_{-}+ \frac{d^2f}{dz_{+}dz_{-}}\mid_{z=Z}B_{+}B_{-}. 
\end{equation}
We can also define the antipode, counite and involution on $C^\infty (R^2)$.
They are given by 
\begin{equation}
\label{27}S(f(z))=f(-z),\ \ \ \varepsilon (f(z))=f(0),\ \ \ (f(z))^{*}=%
\overline{f(z)}, 
\end{equation}
where the bar denotes the usual complex conjugation.

Let $A(E(1,1\mid p))$ be the space of polynomials of $\eta_\pm$ and $\delta$%
. The restrictions (\ref{15}), (\ref{19}) together with (\ref{1}),(\ref{2})
and (\ref{4}) imply that it is finite $*$--Hopf algebra with dimension $p^3$%
. We call it reduced quantum Poincar\'{e} group and denote by $E(1,1\mid p)$.

\begin{tanim}
Quantum Poincar\'e group $E_q(1,1)$ at roots of unity is the $C^{*}$--
algebra $A(E_q(1,1))=A(E(1,1\mid p))\times C^\infty (R^2)$ with a Hopf
algebra structure given by (\ref{2}), (\ref{26}) and (\ref{27}).
\end{tanim}

Let us define the homomorphism $\xi_C : A(E_q(1,1))\rightarrow C^\infty(R^2)$%
, such that 
$$
\xi_C (\eta _{\pm }) = 0,\ \ \ \xi_C (\delta )=1, \ \ \ \xi_C (z_\pm )=
z_\pm . 
$$
From (\ref{26}) we get 
\begin{equation}
\label{tran}\xi_C \diamond \xi_C (f (z) )=f(Z). 
\end{equation}
The operations (\ref{27}) and (\ref{tran}) define a Hopf algebra structure
on $C^\infty(R^2)$. The transformation law 
$$
\xi_C \diamond \xi_C (z_\pm )= z_{\pm }\otimes 1+1\otimes z_{\pm } 
$$
implies that the $*$--Hopf algebra $C^\infty(R^2)$ is the space of all
infinitely differentiable functions on the translation group $R^2$. The
quantum Poincar\'{e} group $E_q(1,1)$ at roots of unity contains the
invariant discrete $E(1,1\mid p)$ and translation $R^2$ subgroups. Using the
group multiplication law (\ref{26}) and analogies with the supersymmetry
theory we call $E_q(1,1)$ $Z_p$--graded product of $E(1,1\mid p)$ and $R^2$.

The quantum group $E(1,1\mid p)$ contains $p$-dimensional invariant subgroup 
$SO(1,1\mid p)$, which is the $*$--Hopf algebra $A(SO(1,1\mid p))$ of
polynomials of $\delta$ subject to the restriction (\ref{15}). The right
sided coset $M^{(1,1)}=E(1,1\mid p) /SO(1,1\mid p)$ is the subspace $%
A(M^{(1,1)})$ of $A(E(1,1\mid p))$ defined as 
$$
A(M^{(1,1)})= \{ a\in A(E(1,1\mid p)): \ \ \  \xi_S\diamond a = a\otimes 1\},
$$
where $\xi_S$ be the homomorphic map of $A(E(1,1\mid p))$ onto $A(SO(1,1\mid
p))$, such that 
$$
\xi_S (\eta _{\pm })=0,\ \ \ \xi_S (\delta )=\delta. 
$$
One can show that 
$$
\xi_S\diamond \eta _{+}^n\eta _{-}^m\delta ^k = \eta _{+}^n\eta
_{-}^m\delta^k\otimes \delta^k%
$$
which implies that  $\eta _{+}^n\eta _{-}^m$, $n$, $m\in [0, p-1]
$, form a basis of $A(E^{(1,1)}_p)$. The elements 
\begin{equation}
e^\pm_{nm}=\frac{\eta^{p-1-n}_+\eta^{p-1-m}_- \pm \eta^n_+\eta^m_- } {\sqrt{%
q^{2n+1}+q^{-2n-1}}}, \ \ \ \ \ n,m\in [0,p-1] 
\end{equation}
also form a basis in $M^{(1,1)})$ which are independent in the range 
$$
n\in [0, n_0-1], \ \ \ m\in [0,2n_0 ] \ \ \ {\rm and } \ \ \ n=n_0, \ \ \
m\in [0,n_0 ], 
$$
where $p=2n_0+1$. The number of independent vectors $e^+_{nm}$ and $e^-_{nm}$
are $\frac{p^2+1}{2}$ and $\frac{p^2-1}{2}$ respectively.

The quantum plane $E^{(1,1)}_q$$=E_q(1,1)/SO(1,1\mid p)$ is the subspace $%
A(E^{(1,1)}_q)$ of $A(E_q(1,1))$ defined as 
$$
A(E^{(1,1)}_q)=A(M^{(1,1)}_p)\times C^\infty (R^2).%
$$

\begin{tanim}
The quantum algebra $U_q(e(1,1))$ at roots of unity is the 
$*$--Hopf algebra generated by $p_{\pm }$ and $\kappa $ 
subject to condition (\ref{15a}). The monomials 
\begin{equation}
\label{32}P_{+}^tP_{-}^sp_{+}^np_{-}^m\kappa ^k,\ \ \ n,\ m,\ k\in [0,\
p-1],\ \ \ t,\ s\in N,
\end{equation}
where $P_{\pm }=p_{\pm }^p$, form a basis of $U_q(e(1,1))$. The $*$--Hopf
algebra structure of $U_q(e(1,1))$ is given by (\ref{5}), (\ref{6}), 
(\ref{8}) and 
$$
\Delta (P_{\pm })=P_{\pm }\otimes 1+1\otimes P_{\pm },\ \ S(P_{\pm
})=-P_{\pm },\ \ \varepsilon (P_{\pm })=0,\ \ P_{\pm }^{*}=P_{\pm }. 
$$
\end{tanim}

The $*$--Hopf algebra $U(r^2)$ generated by $P_{\pm }$ forms the invariant $*
$--sub--Hopf algebra of $U_q(e(1,1))$, which is dual to  the Hopf algebra $%
C^\infty (R^2)$. More precisely due to the Schwartz theorem $U(r^2)$ is
isomorphic to the subspace of distributions on $C^\infty (R^2)$ with support
at the unit element $(0,\ 0)\in R^2$.

The homomorphism $\xi^\prime_C : U_q(e(1,1))\rightarrow U(e(1,1\mid p))$
given by 
$$
\xi^\prime_C (p _\pm ) = p_\pm ,\ \ \ \xi^\prime_C (\kappa )=\kappa , \ \ \
\xi^\prime_C ( P_\pm )= 0 
$$
defines another sub--Hopf algebra of $U_q(e(1,1))$, which is generated by
the elements $p_\pm$ and $\kappa$ subject to the conditions 
$$
p _\pm^p =0,\ \ \ \kappa^p =1_U. 
$$
$U(e(1,1\mid p))$ is in non--degenerate duality with $A(E(1,1\mid p))$.

\vspace{1cm} \noindent 
{\bf 3. Invariant Integral on $E_q(1,1)$}

\vspace{.2cm} \noindent 

\begin{teo}
The linear functional ${\cal I}$ on $A(E(1,1\mid p))$ such that 
$$
{\cal I}(\eta _{+}^n\eta _{-}^m\delta ^k)=q^{-1}\delta _{n,p-1}\delta
_{m,p-1}\delta _{k,0({\rm mod}\ p)} 
$$
defines the unique invariant integral on the reduced quantum Poincar\`e
group $E(1,1\mid p)$.
\end{teo}
$Proof$. Let us find the linear functional ${\cal I}^{\prime }$ on $%
A(E(1,1\mid p))$ which for any element $a$ from $A(E(1,1\mid p))$ satisfies
the left 
$$
{\cal I}^{\prime }\diamond a={\cal I}^{\prime }(a)1_A 
$$
and right 
$$
a\diamond {\cal I}^{\prime }={\cal I}^{\prime }(a)1_A 
$$
invariance conditions. By the virtue of (\ref{13}) for $a=\eta _{+}^n\eta
_{-}^m\delta ^k$ the left invariance condition reads 
$$
\sum_{t,s=0}^{n,m}\left[ 
\begin{array}{c}
n \\ 
t
\end{array}
\right] _{+}\left[ 
\begin{array}{c}
m \\ 
s
\end{array}
\right] _{-}q^{2t(s-m)}\eta _{+}^{n-t}\eta _{-}^{m-s}\delta ^{k-s+t}{\cal I}%
^{\prime }(\eta _{+}^t\eta _{-}^s\delta ^k)=1_A{\cal I}^{\prime }(\eta
_{+}^n\eta _{-}^m\delta ^k) 
$$
which implies 
$$
{\cal I}^{\prime }\diamond (\eta _{+}^t\eta _{-}^s\delta ^k)=0\ \ \ {\rm for}%
\ \ \ t\in [0,n-1],\ s\in [0,m-1] 
$$
and 
\begin{equation}
\label{i1}k+n-m=0({\rm mod}\ p).
\end{equation}
If $n,\ m\in [0,p-2]$ we can employ the above reasoning for the element $%
a=\eta _{+}^{n+1}\eta _{-}^{m+1}\delta ^k$ and obtain 
\begin{equation}
\label{i2}{\cal I}^{\prime }(\eta _{+}^n\eta _{-}^m\delta ^k)=0\ \ \ {\rm for%
}\ \ \ n,\ m\in [0,p-2].
\end{equation}
(\ref{i1}) and (\ref{i2}) imply that the linear functional ${\cal I}^{\prime
}$ satisfies the left invariance condition if 
$$
{\cal I}^{\prime }(\eta _{+}^n\eta _{-}^m\delta ^k)=\omega \delta
_{n,p-1}\delta _{m,p-1} \delta_{k,0({\rm mod}\ p)}, 
$$
where $\omega $ is an arbitrary complex number. In a similar fashion one can
show that the right invariance implies the same condition on ${\cal I}%
^{\prime }$. Thus every linear functional on $A(E(1,1\mid p))$ satisfying
the left and right invariance conditions is proportional to ${\cal I}$. 
\hspace{.5cm} $\Box $

\vspace{.2cm} 
\noindent
Define the bilinear form $(\cdot ,\cdot )_p$ on $E(1,1\mid p)$ by 
\begin{equation}
\label{I17}(a,b) = {\cal I} (ab^* ). 
\end{equation}
Because of  the property 
$$
{\cal I} ( a^* )= \overline{{\cal I} (a )}%
$$
this bilinear form is Hermitian. The vectors $e^\pm_{nm}$ spanning the basis
of the coset space $A(M^{(1,1)})$ are orthonormal with respect to the above
form 
\begin{equation}
\label{ort}(e^\pm_{nm},e^\pm_{n^\prime m^\prime})=\pm
\delta_{nn^\prime}\delta_{mm^\prime}, \ \ \ (e^\pm_{nm},e^\mp_{n^\prime
m^\prime})=0. 
\end{equation}
Thus $A(M^{(1,1)}_p)$ equipped with the Hermitian form (\ref{I17}) is the
pseudo--Euclidean space with $\frac{p^2+1}{2}$ positive and $\frac{p^2-1}{2}$
negative signatures.

Let ${\cal I}_C$ be the linear functional on the space $C^\infty(R^2)$ of
all infinitely differentiable functions with finite support in $R^2$ given
by 
\begin{equation}
{\cal I}_C (f)= \int_{-\infty}^\infty \int_{-\infty}^\infty dz_+ dz_- f
(z_+,z_- ) 
\end{equation}
and let $A_0(E_q(1,1))$ be the subspaces 
$$
C_0^\infty(R^2)\times A(E(1,1\mid p))%
$$
of $A(E_q(1,1))$ whose any element $F$ is the finite sum 
$$
F=\sum_na_nf_n, 
$$
where $f_n\in C_0^\infty (R^2)$ and $a_n\in A(E(1,1\mid p))$. It is clear
that ${\cal I}_C$ is the invariant integral on the translation group
satisfying the properties 
\begin{equation}
\label{inv}({\cal I}_C\otimes id)(\xi_C\diamond\xi_C)(f)= \zeta (f), \ \ \
(id \otimes {\cal I}_C)(\xi_C\diamond\xi_C)(f)= \zeta (f) 
\end{equation}
for any $f\in C_0^\infty(R^2)$.

\begin{teo}
The linear functional ${\cal I}_E$ on $A_0(E_q(1,1))$ given by 
$$
{\cal I}_E(F)=\sum_n{\cal I}(a_n){\cal I}_C(f_n) 
$$
defines the unique invariant integral on the quantum Poincar\`e group $%
E_q(1,1)$.
\end{teo}
$Proof.$ By the virtue of (\ref{26}) and (\ref{tran}) for $G=af$ we have 
\begin{eqnarray*}
{\cal I}_E\diamond G  & = & (id\otimes {\cal I}_E)  [ \Delta (a) 
\{ (\xi_C\diamond\xi_C)(f)+ B_+(\xi_C\diamond\xi_C)(\frac{df}{dz_+}) \\
& & + B_- (\xi_C\diamond\xi_C)(\frac{df}{dz_-}) + B_+B_- 
(\xi_C\diamond\xi_C)(\frac{d^2f}{dz_+dz_-}) \}].
\end{eqnarray*}
By making use of (\ref{inv}) we get 
\begin{eqnarray*}
{\cal I}_E\diamond G & = & 1_A {\cal I}(a){\cal I}_C (f)+
(id\otimes {\cal I})[ \Delta (a)
\{ B_+ {\cal I}_C(\frac{df}{dz_+}) + 
B_-  {\cal I}_C (\frac{df}{dz_-})\}] \\
& &+  (id\otimes {\cal I}) (\Delta (a)B_+B_- )
{\cal I}_C (\frac{d^2f}{dz_+dz_-}).
\end{eqnarray*}
Using the properties 
$$
\zeta _C(\frac{df}{dz_{\pm }})=0,\ \ \ \zeta _C(\frac{d^2f}{dz_{+}dz_{-}})=0 
$$
satisfied by the functions $f\in C_0^\infty (R^2)$ we arrive at 
$$
{\cal I}_E\diamond G=1_A{\cal I}(a){\cal I}_C(f)=1_A{\cal I}_E(G), 
$$
which together with the linearity of the functional ${\cal I}_E$ implies 
$$
{\cal I}_E\diamond F=1_A{\cal I}_E(F) 
$$
for any $F\in A_0(E_q(1,1))$. We have proved the left invariance condition.
In a similar fashion one can prove the right invariance condition. The
uniqueness of the invariant integral ${\cal I}_E$ follows from the
uniqueness of the invariant integrals ${\cal I}$ and ${\cal I}_C$. 
\hspace{.5cm} $\Box $

\vspace{.2cm} \noindent
By means of the invariant integral we define in $E_q(1,1)$ the bilinear 
form by 
\begin{equation}\label{her}
(F,G)_E={\cal I}_E(FG^{*}), 
\end{equation}
where $F$, $G\in A_0(E_q(1,1))$. Because of the property 
$$
{\cal I}_E(F^{*})=\overline{{\cal I}_E(F)}%
$$
this bilinear form is Hermitian.

Let $A_0(E_q^{(1,1)})$ be the subspace 
$$
C_0^\infty(R^2)\times A(M^{(1,1)}).%
$$
of $A(E_q^{(1,1)})$ whose any element $X$ is the finite sum 
$$
X=\sum_{nm}f_{nm}^{+}e_{nm}^{+}+ \sum_{nm}f_{nm}^{-}e_{nm}^{-}, 
$$
where $e_{nm}^{\pm}$ form a basis of $A(M^{(1,1)})$ and $f_{nm} \in
C_0^\infty(R^2)$. By the virtue of (\ref{ort}) we get 
\begin{equation}
(X,X)_E=\sum_{nm} {\cal I}_C(f_{nm}^{+}\overline{f_{nm}^{+}})- \sum_{nm}%
{\cal I}_C(f_{nm}^{-}\overline{f_{nm}^{-}}), 
\end{equation}
which implies that $A_0(E_q^{(1,1)})$ equipped with the Hermitian form (\ref
{her}) is the pseudo--Euclidean space.

\vspace{1cm} \noindent 
{\bf 4. Irreducible $*$--Representations of $U_q(e(1,1))$ }

\vspace{.2cm} \noindent 
The homomorphism ${\cal L}^\lambda $ : $U_q(e(1,1))\rightarrow $ Lin $%
A(SO(1,1\mid p))$ given by 
\begin{equation}
\label{rep}{\cal L}^\lambda (p_{\pm })\delta ^m=\lambda _{\pm }\delta ^{m\pm
1},\ \ \ {\cal L}^\lambda (\kappa )\delta ^m=q^m\delta ^m
\end{equation}
for $\lambda =(\lambda _{+},\ \lambda _{-})\neq (0,\ 0)$ defines
p--dimensional irreducible representation of the quantum algebra $U_q(e(1,1))
$ in the linear space $A(SO(1,1\mid p))$. Since $\delta ^p=1_A$ for any $%
a\in A(SO(1,1\mid p))$ we have 
$$
{\cal L}^\lambda (P_{\pm })a=\lambda _{\pm }^pa 
$$
This representation is cyclic. For $\lambda =(0,\ 0)$ we have one
dimensional representation 
\begin{equation}
{\cal L}^{(m)}(p_{\pm })\delta ^m=0,\ \ \ \ {\cal L}^{(m)}
(\kappa )\delta^m=q^m\delta ^m
\end{equation}
with the weight $m\in [0,\ p-1]$. The homomorphisms ${\cal L}^\lambda $ and $%
{\cal L}^{(m)}$ exhaust all irreducible representations of the quantum 
algebra $U_q(e(1,1))$. This is rather trivial consequence of the general theory
presented in \cite{kac}, to which we refer for proof and details. 
Representations of the quantum algebra $U_q(e(1,1))$ is also considered in 
\cite{cic}. However the quantum algebra studied in \cite{cic} differs 
because there the restriction (\ref{15}) is not considered.

Let us find out when the homomorphism ${\cal L}^\lambda$ defines $*$%
--representation of the quantum algebra $U_q(e(1,1))$, that is when for any $%
\phi \in U_q(e(1,1))$ we have 
\begin{equation}
({\cal L}^\lambda (\phi ) )^* = {\cal L}^\lambda (\phi^*) 
\end{equation}
For this purpose we define in $A(SO(1,1\mid p))$ the Hermitian form 
\begin{equation}
\label{hers}(a,b)_S={\cal I}_S (a^* b ), 
\end{equation}
where ${\cal I}_S$ is the invariant integral on $SO(1,1\mid p)$ given by 
$$
{\cal I}_S ( \delta^m )=\delta_{m, 0({\rm mod} \ p)}.
$$
For $n$, $m\in [0, p-1]$ we have 
\begin{equation}
\label{a}(\delta^n, \delta^m)_S= \delta_{m+n, 0}+ \delta_{m+n, p}, 
\end{equation}
which implies that the vectors 
$$
e_m^{\pm }=\frac {1}{\sqrt{2}} (\delta^m \pm \delta^{p-m} ), \ \ \ \ \ m\in[%
0,\frac{p-1}{2}] 
$$
are orthonormal with respect to the Hermitian form (\ref{hers}) 
$$
(e_m^\pm ,e_k^\pm )_S= \pm \delta _{mk},\ \ \ \ (e_m^\mp ,e_k^\pm )_S=0. 
$$
The $*$--Hopf algebra $A(SO(1,1\mid p))$ equipped with the Hermitian form (%
\ref{hers}) is pseudo--Euclidean space with $\frac{p+1}{2}$ positive and $%
\frac{p-1}{2}$ negative signatures.

The adjoint $({\cal L}^\lambda (\phi ))^*$ of the linear operator ${\cal L}%
^\lambda (\phi )$ is defined as 
$$
({\cal L}^\lambda (\phi )a,b)_S= (a, ({\cal L}^\lambda (\phi ) )^* b)_S, 
$$
where $a, \ b$ are arbitrary elements from $A(SO(1,1\mid p))$. Using the
representation formula (\ref{rep}) and the involution (\ref{8}) we conclude
that when $\lambda_\pm $ are real numbers the homomorphism ${\cal L}^\lambda$
defines $*$--representation of the quantum algebra $U_q(e(1,1))$. The
homomorphism ${\cal L}^{(m)}$ also defines $*$--representation of $U_q(e(1,1))$.

\vspace{1cm} \noindent 
{\bf 5. Pseudo--Unitary Irreducible Representations of $E_q(1,1)$ }

\vspace{.2cm} \noindent 
Let us briefly recall the construction and the main properties of  $%
universal\ T--matrix$ \cite{fron}. Consider two Hopf algebras $A(G)$ and $%
U(g)$ in non--degenerate duality. Let $\{x_a\}$ and $\{X^b\}$ be dual linear
basis of $A(G)$ and $U(g)$ respectively, with $a$ and $b$ running in an
appropriate set of indices, so that $\langle x_a,X^b\rangle =\delta _{ab}$.
We define the element $T\in U(g)\otimes A(G)$ as 
$$
T=\sum_ax_a\otimes X^a. 
$$
The universal $T$--matrix is a resolution of the identity which maps the Lie
group $G$ into itself. Moreover, if we choose the representation of $U(g)$
we correspondingly obtain the corepresentation of $A(G)$ or representation
of $G$.

The elements $z_{+}^tz_{-}^s\eta _{+}^n\eta _{-}^m\zeta (k)$ and (\ref{32})
defines the linear basis in $A(E_q(1,1))$ and $U_q(e(1,1))$ respectively.
Introducing the cut off q--exponential 
\begin{equation}
\label{40}e_{\pm }^x=\sum_{m=0}^{p-1}\frac{q^{\pm \frac{m(m-1)}2}}{[m]!}x^m.
\end{equation}
by the direct calculation we arrive at the following result.

\begin{pro}
We have the duality relations 
\begin{eqnarray*}
\langle P_+^t P_-^s p_+^n p_-^m \kappa^k, 
z_+^{t^\prime }z_-^{s^\prime } \eta_+^{n^\prime } 
\eta_-^{m^\prime} \zeta (k^\prime )\rangle   
& = & i^{n+m+t+l} q^{\frac{n-m}{2}-nm}t!s![n]![m]! \\
& & \delta_{nn^\prime } \delta_{mm^\prime }
\delta_{tt^\prime } \delta_{ll^\prime }
\delta_{k+t+l, k^\prime },  
\end{eqnarray*}
which implies that the universal $T$--matrix in $U_q(e(1,1))\otimes
A(E_q(1,1))$ has the form 
$$
T=e^{-iP_{+}\otimes z_{+}-iP_{-}\otimes z_{-}}e_{+}^{i\epsilon _{+}\otimes
\eta _{+}}e_{-}^{i\epsilon _{-}\otimes \eta _{-}}D(\kappa ,\delta ), 
$$
where 
$$
\epsilon _{\pm }=-q^{\mp \frac 12}p_{\pm }\kappa ^{-1} 
$$
and 
$$
D(\kappa ,\delta )=\frac 1p\sum_{m,k=0}^{p-1}q^{-mk}\kappa ^m\otimes \delta
^k 
$$
\end{pro}
The universal $T$--matrix satisfies the properties 
\begin{equation}
\label{41}[(*\otimes *)T]\cdot T=1_U\otimes 1_A,\ \ \ T\cdot
[(*\otimes*)T]=1_U\otimes 1_A 
\end{equation}
and 
\begin{equation}
\label{42}(id\otimes \Delta )T=(T\otimes 1_A)(id\otimes \sigma )(T\otimes
1_A), 
\end{equation}
where $\sigma (F\otimes G)=G\otimes F$, $F$, $G\in A(E_q(1,1))$ is the
permutation operator.

Define the linear map $T^\lambda : A(SO(1,1\mid p)) \rightarrow A(SO(1,1\mid
p))\otimes A(E_q(1,1))$, such that 
\begin{equation}
\label{repg}T^\lambda a =e^{-i{\cal L}^\lambda (P_+)\otimes z_+ - i{\cal L}%
^\lambda (P_-)\otimes z_-} e_+^{i {\cal L}^\lambda (\epsilon_+)\otimes
\eta_+ } e_-^{i {\cal L}^\lambda (\epsilon_-)\otimes \eta_- } D({\cal L}%
^\lambda(\kappa ) , \delta ) (a\otimes 1). 
\end{equation}
Due to (\ref{42}) and the irreducibility of the representation ${\cal L}%
^\lambda$ we conclude that the above linear map defines $p$--dimensional
irreducible representations of the quantum Poincar\' {e} group in the linear
space $A(SO(1,1\mid p))$. Let us extend the Hermitian form (\ref{hers}) to
the form $\{ \cdot , \cdot \}_S$ by setting 
\begin{equation}
\label{herss}\{ a\otimes F, b\otimes G \}_S = F^*G(a,b)_S, 
\end{equation}
where $F$, $G\in A(E_q(1,1))$ and $a$, $b\in A(SO(1,1\mid p))$. When $%
\lambda_\pm$ are real numbers due to (\ref{41}) we get 
\begin{equation}
\label{pseudo}\{ T^\lambda a ,T^\lambda b \}_S = (a,b)_S 1_A. 
\end{equation}
Thus the irreducible representation $T^\lambda$ of the quantum group $%
E_q(1,1)$ in the pseudo--Euclidean space $A(SO(1,1\mid p))$ is
pseudo--unitary when $\lambda_\pm\in R$.

By the virtue of the representation formula (\ref{repg}) and the relation (%
\ref{a}) we obtain the integral representation for the matrix elements of
the irreducible pseudo--unitary representations $T^\lambda $ 
\begin{equation}
D_{mn}^\lambda =\{\delta ^{p-m}\otimes 1_A,T^\lambda \delta ^n\}_S.
\end{equation}
After lengthily but straightforward calculations we have the following
result.

\begin{pro}
The matrix elements of the pseudo--unitary irreducible representations of $%
E_q(1,1)$ are 
$$
D_{mn}^\lambda =e^{-i\lambda _{+}^pz_{+}-i\lambda
_{-}^pz_{-}}[\sum_{k=0}^{p-1-n+m}\frac{(-\lambda ^2)^kq^{-k(m+n)}}{%
[k]![k+n-m]!}\xi ^k(-iq^{(\frac 12-n)}\lambda _{-}\eta _{-})^{n-m}\delta ^n 
$$
$$
+(-iq^{(-\frac 12-n)}\lambda _{+}\eta _{+})^{p+m-n}\delta ^n\sum_{k=0}^{n-m}%
\frac{(-\lambda ^2)^kq^{k(m+n)}}{[k]![k+p+m-n]!}\xi ^k]\ \ \ for\ \ \ n\geq
m 
$$
and 
$$
D_{mn}^\lambda =e^{-i\lambda _{+}^pz_{+}-i\lambda
_{-}^pz_{-}}[\sum_{k=0}^{m-n}\frac{(-\lambda ^2)^kq^{-k(m+n)}}{[k]![k+p+n-m]!
}\xi ^k(-iq^{(\frac 12-n)}\lambda _{-}\eta _{-})^{p+n-m}\delta ^n 
$$
$$
+(-iq^{(-\frac 12-n)}\lambda _{+}\eta _{+})^{m-n}\delta
^n\sum_{k=0}^{p-1-m+n}\frac{(-\lambda ^2)^kq^{k(m+n)}}{[k]![k+m-n]!}\xi ^k]\
\ \ for\ \ \ m\geq n, 
$$
where $\xi =q\eta _{+}\eta _{-}$ and $\lambda ^2=\lambda _{+}\lambda _{-}$.
\end{pro}

For the special case $D_{m0}^\lambda $ we have the explicit formula 
\begin{equation}
\label{mat}D_{m0}^\lambda =e^{-i\lambda _{+}^pz_{+}-i\lambda _{-}^pz_{-}}[%
{\cal J}_{p-m}(\lambda ^2\xi )(-iq^{\frac 12}\lambda _{-}\eta
_{-})^{p-m}+(-iq^{-\frac 12}\lambda _{+}\eta _{+})^m{\cal J}_m(\lambda ^2\xi
)],
\end{equation}
where $m\in [0,\ p-1]$ and 
\begin{equation}
\label{bes}{\cal J}_m(x)=\sum_{k=0}^{p-1-m}\frac{(-1)^k}{[k]![k+m]!}(q^mx)^k.
\end{equation}
The pseudo--unitarity condition (\ref{pseudo}) implies 
\begin{equation}
(D_{0m}^\lambda )^{*}D_{0n}^\lambda +\sum_{k=1}^{p-1}(D_{km}^\lambda
)^{*}D_{p-kn}^\lambda =(\delta ^m,\delta ^n)_S1_A.
\end{equation}
Special cases are 
$$
(D_{00}^\lambda )^{*}D_{00}^\lambda +\sum_{k=1}^{p-1}(D_{k0}^\lambda
)^{*}D_{p-k0}^\lambda =1_A 
$$
and 
$$
(D_{0s}^\lambda )^{*}D_{0p-s}^\lambda +\sum_{k=1}^{p-1}(D_{ks}^\lambda
)^{*}D_{p-kp-s}^\lambda =1_A, 
$$
where $s\in [1,p-1]$. Moreover, we have the addition theorem 
\begin{equation}
\Delta (D_{nm}^\lambda )=\sum_{k=0}^{p-1}D_{nk}^\lambda \otimes
D_{km}^\lambda .
\end{equation}
The pseudo--unitary representation $T^{(m)}$ of the quantum Poincar\`e group
corresponding to the $*$--representation ${\cal L}^m$ is given by 
$$
T^{(m)}\delta ^m=\delta ^m\otimes \delta ^m, 
$$
where $m\in [0,\ p-1]$.

\vspace{.2cm} \noindent
$Remarks$. (i) Recall that the Hahn--Exton q--Bessel functions $J_m(x)$
related to the unitary irreducible representations of the quantum Euclidean
group $E_q(2)$ are \cite{koel} 
$$
J_m(x)=\sum_{k=0}^\infty\frac {(-1)^k}{[k]![k+m]!}(q^mx)^k.%
$$
Comparing (\ref{bes}) to the above expression we conclude that the matrix
elements of the pseudo--unitary irreducible representations of the quantum
Poincar\'e group are the cut off Hahn--Exton q--Bessel function.

(ii) Inspecting (\ref{repg}) we observe that irreducible representations of $%
E_q(1,1)$ are induced by the irreducible representations of the translation 
subgroup $R^2$.

(iii) The linear map $T^{(m)}$ defines the one dimensional pseudo--unitary
representations of the invariant subgroup $SO(1,1\mid p)\in E_q(1,1)$.

\vspace{1cm} \noindent 
{\bf 6. Quasi--Regular Representation}

\vspace{.2cm} \noindent 
The comultiplication 
\begin{equation}
\label{quasi}\Delta : A_0(E_q^{(1,1)})\rightarrow A_0(E_q(1,1))\otimes
A_0(E_q^{(1,1)}) 
\end{equation}
defines the left quasi--regular representation of the quantum Poincar\'e
group $E_q(1,1)$ in the vector space $A_0(E_q^{(1,1)})$. Let us extend the
Hermitian form (\ref{her}) to the form $\{\cdot , \cdot\}_E$ by setting 
$$
\{ F\otimes X , G\otimes Y\}_E = FG^* (X, Y)_E, 
$$
where $X$, $Y\in A_0(E_q^{(1,1)})$ and $F$, $G\in A_0(E_q(1,1))$. Since the
Hermitian form $(\cdot , \cdot )_E$ is defined by means of the invariant
integral we have 
\begin{equation}\label{here}
\{ \Delta (X) , \Delta (Y)\}_E =1_A (X, Y)_E, 
\end{equation}
which implies that the left quasi--regular representation (\ref{quasi}) is
pseudo--unitary.

The right representation ${\cal R}$ of the quantum algebra $U_q(e(1,1))$
corresponding to the left quasi--regular representation (\ref{quasi}) is
given by 
$$
{\cal R}(\phi ) F = F\diamond \phi . 
$$
We have 
\begin{equation}
\label{dif1}{\cal R}( p_\pm ) \eta_\pm^k =i q^{\pm \frac{1}{2}}[k]
\eta_\pm^{k-1}, \ \ \ {\cal R}( p_\pm ) \eta_\mp^k =0, \ \ \ {\cal R}(
\kappa )\eta_\pm^k = q^{\pm k}\eta_\pm^k 
\end{equation}
and 
\begin{equation}
\label{dif2}{\cal R}( p_\pm ) f= i q^{\pm \frac{1}{2}} \frac{(-1)^{\frac{p+1%
}{2}} }{[p-1]!} \eta^{p-1}_\pm \frac{df}{dz_\pm }, \ \ \ {\cal R}( P_\pm )
f= i \frac{df}{dz_\pm },  \ \ \ {\cal R}(\kappa ) f= f, 
\end{equation}
where $f\in C_0^\infty (R^2)$. Using the following relations satisfied by
the right representation ${\cal R}$ 
$$
{\cal R}( \phi \phi^\prime )={\cal R}(\phi^\prime ) {\cal R}(\phi ) , 
$$
$$
{\cal R}( p_\pm ) (X Y )= {\cal R}( p_\pm )X {\cal R}(\kappa ) Y + {\cal R}%
(\kappa^{-1})X {\cal R}(p_\pm )Y, 
$$
$$
{\cal R}(\kappa ) (XY)= {\cal R}(\kappa )X {\cal R}(\kappa )Y 
$$
we can define the action of an arbitrary operator ${\cal R}(\phi )$ on any
function from $A_0(E_q^{(1,1)})$. Due to the identity 
$$
\overline{\langle \phi , F^*\rangle }= \langle (S(\phi))^* ,F\rangle, \ \ \
\ \ F\in A_0(E_q(1,1)) 
$$
and the pseudo--unitarity condition (\ref{here}) for any $\phi\in U_q(e(1,1))
$ we have 
$$
({\cal R}(\phi )X, Y )_E= (X, {\cal R}(\phi^*)Y )_E 
$$
Thus the antihomomorphism ${\cal R}$ : $U_q(e(1,1))\rightarrow$ Lin $%
A_0(E_q^{(1,1)})$ defines $*$--representation of the quantum algebra $%
U_q(e(1,1))$ in the pseudo--Euclidean space $A_0(E_q^{(1,1)})$.

The quantum algebra $U_q(e(1,1))$ has three Casimir elements $P_\pm$ and $%
p_+p_-$ with one restriction 
$$
P_+P_- = (p_-p_+)^p. 
$$
Therefore irreducible representations of $U_q(e(1,1))$ will be labelled by
two indices. We construct the irreducible representations of the quantum
algebra $U_q(e(1,1))$ in the pseudo--Euclidean space $A_0(E_q^{(1,1)})$ by
diagonalizing the complete set of commuting elements of $U_q(e(1,1))$ in $%
A_0(E_q^{(1,1)})$.

\vspace{1cm} \noindent
(i) $The \ angular \ momentum \ states$. Choose the following complete set
of observables : ${\cal R}(P_{\pm })$, ${\cal R}(p_{+}p_{-})$, ${\cal R}%
(\kappa )$. Inspecting (\ref{dif1}) and (\ref{dif2}) we observe that the
functions 
$$
X=e^{-i\lambda _{+}^pz_{+}-i\lambda _{-}^pz_{-}} [X_1(\xi )\eta
_{-}^{p-m}+\eta_{+}^m X_2(\xi )], 
$$
with $X_1(\xi )$ and $X_2(\xi )$ being some polynomials, are 
eigenstates of the linear operators ${\cal R}(P_{\pm })$ and ${\cal R}%
(\kappa )$ with eigenvalues $\lambda _{\pm }^p$ and $q^m$ respectively. 
The eigenvalue equation 
$$
{\cal R}(p_{+}p_{-})X=\lambda ^2X%
$$
is solved by 
$$
X=D_{m0}^\lambda ,
$$
where $\lambda ^2=\lambda _{+}\lambda _{-}$ and $D_{m0}^\lambda $ are the
matrix elements (\ref{mat}). By direct calculations we arrive at
the following results.

\begin{pro}
The right representation of $U_q(e(1,1))$ on the matrix elements $%
D_{m0}^\lambda $ is given by 
$$
{\cal R}(p_{+})D_{m0}^\lambda =\lambda _{+}D_{m-1,0}^\lambda ,\ \ \ \ \ m\in
[1,p-1], 
$$
$$
{\cal R}(p_{-})D_{m0}^\lambda =\lambda _{-}D_{m+1,0}^\lambda ,\ \ \ \ \ m\in
[0,p-2] 
$$
and 
$$
{\cal R}(p_{+})D_{00}^\lambda =\lambda _{+}D_{p-10}^\lambda ,\ \ \ \ {\cal R}%
(p_{-})D_{p-1,0}^\lambda =\lambda _{-}D_{00}^\lambda . 
$$
\end{pro}

\begin{pro}
The matrix elements of the irreducible pseudo--unitary representation
satisfy the orthogonality condition 
$$
(D_{n0}^\lambda ,D_{m0}^{\lambda ^{\prime }})_E=\Lambda
\delta (\lambda _{+}-\lambda
_{+}^{\prime })\delta (\lambda _{-}-\lambda _{-}^{\prime })
\delta _{n+m,0({\rm mod} \ p)}, 
$$
where 
$$
\Lambda=\frac{2\pi}{p^2}\sum_{k=0}^{p-1}\frac 1{([k]![p-1-k]!)^2}. 
$$
is the normalization constant. 
\end{pro}

\vspace{.5cm} \noindent
(ii) $The \ Plane \ wave \ states$. We choose the following complete set of
observables: ${\cal R}(P_\pm )$, ${\cal R}(p^\prime_\pm )$, where 
$$
p_+^\prime =q^{-\frac{1}{2}} p_+\kappa^{-1}, \ \ \ \ p_-^\prime = q^{-\frac{1%
}{2}} p_-\kappa. 
$$
Due to the relation $P_\pm =-(p_\pm^\prime )^p$ it is sufficient to solve
the eigenvalue equations 
\begin{equation}
\label{ev}{\cal R}(p_\pm^\prime )Y= \chi_\pm Y. 
\end{equation}
\begin{pro}
The eigenfunctions of (\ref{ev}) are 
$$
Y=e_{+}^{-i\chi _{+}\eta _{+}}e_{+}^{-iq\chi _{-}\eta _{-}}e^{i\chi
_{+}^pz_{+}}e^{i\chi _{-}^pz_{-}}, 
$$
where $e_{+}^x$ is the cut off exponential (\ref{40}).
\end{pro}
$Proof.$ Substituting 
$$
Y=e^{i \chi^p_+ z_+}e^{i\chi^p_- z_-}Y_+(\eta_+) Y_-(\eta_-)%
$$
in (\ref{ev}) we get 
$$
[{\cal R}(p_+^\prime ) - q\frac{(-1)^{\frac{p+1}{2}} }{[p-1]!}
\chi_+^p\eta^{p-1}_+ ]Y_+ = \chi_+ Y_+ 
$$
and 
$$
[{\cal R}(p_-^\prime ) - \frac{(-1)^{\frac{p+1}{2}} }{[p-1]!}
\chi_-^p\eta^{p-1}_- ]Y_- = \chi_- Y_- ,
$$
which imply the desired result. \hspace{.5cm} $\Box$

\vspace{1cm} \noindent
$Acknowledgement.$ The author thanks Dayi, \"{O}. F. and Duru, I. H. for
helpful discussions.

\end{document}